%Paper: alg-geom/9412012
%From: jml@shire.math.columbia.edu (J. M. Landsberg)
%Date: Tue, 13 Dec 94 16:11:02 -0500
%Date (revised): Fri, 20 Oct 95 10:21:07 EDT

%this is an AMSTeX file
%The maps in the diagram (7.5) are named respectively, \pi and \rho
%The maps in the diagram (7.6) are named
% respectively, \pi_{\Delta} and \rho_{\Delta}

\input amstex
\def\aa#1{A_{#1}}

\def\bcc#1{\Bbb C^{#1}}

\def\cf{\Cal F}
\def\ci{\Cal I}

\def\frp#1#2{\frac{\partial {#1}}{\partial {#2}}}
\def\gg#1#2{g^{#1}_{#2}}

\def\hd{, \hdots ,}

\def\inv{{}^{-1}}
\def\ii{|II|}
\def\hii{\widehat{|II|}}

\def\na{n+a}
\def\ooo#1#2{\omega^{#1}_{#2}}
\def\oo#1{\omega^{#1}}
\def\ot{\!\otimes\!}
\def\pp#1{\Bbb P^{#1}}

\def\ppp{\Bbb P}
\def\pii{\partial II}
\def\qq#1#2#3{q^{#1}_{{#2} {#3}}}
\def\qqn#1#2#3{q^{n+ #1}_{{#2} {#3}}}
\def\ra{\rightarrow}
\def\rr#1#2#3#4{r^{#1}_{{#2} {#3}{#4}}}
\def\rrn#1#2#3#4{r^{n+ #1}_{{#2} {#3}{#4}}}

\def\tann{\text{Ann}}
\def\tsingloc{\text{singloc}}
\def\tbaseloc{\text{baseloc}}
\def\tdim{\text{dim}\,}
\def\tker{\text{ker}\, }

\def\tmod{\text{ mod }}

\def\tm{\tau (X)}\def\sm{\sigma (X)}
\def\tx{\tau (X)}\def\sx{\sigma (X)}
\def\trank{\text{rank }}

\def\up#1{{}^{({#1})}}
\def\upperp{{}^{\perp}}

\def\ww{\wedge}

\define\intprod{\mathbin{\hbox{\vrule height .5pt width 3.5pt depth 0pt %
        \vrule height 6pt width .5pt depth 0pt}}}
\documentstyle{amsppt}
\magnification = 1200
\hsize =14truecm
\hcorrection{.5truein}
\baselineskip =18truept
\vsize =22truecm
\NoBlackBoxes

\topmatter
\title
On Degenerate Secant and Tangential  Varieties and Local Differential Geometry
\endtitle
\rightheadtext{ Degenerate Secant Varieties}
\author
J.M. Landsberg
\endauthor

\address{ Dept. of Mathematics, Columbia University, New York, NY 10027}
\endaddress
\email {jml\@math.columbia.edu }
\endemail
\date {September 4, 1995 }\enddate
\thanks {This work was done while the author was partially supported by an NSF
postdoctoral fellowship and later by  NSF grant DMS-9303704}
\endthanks

\abstract{ We study  the local differential geometry of  varieties
$X^n\subset \Bbb C\Bbb P^{n+a}$ with degenerate secant and tangential
varieties.
We show that the second fundamental form of a smooth variety with degenerate
tangential variety is subject to certain rank restrictions. The
rank restrictions imply
a slightly refined version of Zak's theorem on linear normality,
a short proof of the Zak-Fantecchi
theorem on the superadditivity of higher secant
defects, and a short proof and
extension of a result of Robert's on the nondegeneracy
of secant varieties of Veronese re-embeddings.
We study the geometric structure of   systems
of quadrics generated by the second fundamental
forms of  varieties with degenerate tangential varieties.
In particular, Clifford algebras make an appearance.
 We give a new  proof of Zak's theorem on Severi varieties
that is substantially shorter than the original
by utilizing the Clifford algebra structure and
the rank restrictions.
  We also prove local and global results on the dimensions
of   Gauss images of degenerate tangential varieties. }
\endabstract

\endtopmatter

\document

\noindent\S 1. {\bf Introduction,  Conventions}

One way to study  geometric properties of a variety
$X^n\subset\Bbb C\pp{n+a}$ is by studying coarse geometric
properties of auxiliary varieties one constructs from $X$. The auxiliary
varieties we will study in this paper are the secant variety
$\sx$ and the tangential variety $\tx$, and the coarse properties
of $\sx$ and $\tx$ we will study are their dimensions. For
information on how this study fits into larger questions,
see [LV].  It turns out that
smooth varieties of small codimension with degenerate
$\tx$ (degenerate meaning that $\tx$ is not the entire ambient space)
carry a remarkable amount of infinitesimal geometric structure.
Before going into details, we will need a few definitions.

\smallpagebreak

   Given a variety
$X^n\subset \Bbb P^{{n+a}}$, the {\it secant variety } $\sigma (X)$ of $X$ is
defined to be the union of all points on all secant and tangent lines (i.e.
$\Bbb P^1$'s) of $X$. More precisely,
given  $p,q\in  \Bbb P^{n+a}$,
let $\Bbb P^1_{pq}\subset \Bbb P^{n+a}$ denote the
projective
line  containing $p$ and $q$. Then
$$
\sigma (X) := \overline{\{ x\in \Bbb P^{{n+a}} | x\in \Bbb P^1_{pq}\  \text{for
some} \  p,q \in X\}}.
$$

Secant varieties
have been studied extensively. Two important results on them are:

\proclaim{Zak's Theorem on Linear Normality 1.1} ([FL],[Z1]) Let $X^n\subset
\Bbb C\pp{n+a}$ be a smooth  variety not contained in a hyperplane with $\sigma
(X)\neq\Bbb C \pp{n+a}$. Then $a\geq \frac{n}{2}+2$.

\endproclaim

\proclaim{Zak's Theorem on Severi Varieties 1.2} ([LV],[Z1]) Let   $X^n\subset
\Bbb C\pp{n+a}$
 be a smooth variety  not contained in a hyperplane  with
$\sm \neq
\Bbb C\pp{n+a}$. If $a=\frac{n}{2}+2$, then $X$ is one of
$$\align & i.\ \text{Veronese}\ \pp 2 \subset \pp 5\\
& ii.\ \text{Segre}\  \pp 2\times\pp 2 \subset \pp 8\\
& iii.\ \text{Pl\"ucker  embedded Grassmannian}\ G(\Bbb C^2,\Bbb C^6)\subset
\pp {14}\\
& iv.\  E_6/P \subset \pp {26}.\endalign
$$
\endproclaim

These  four varieties, now called
{\it Severi varieties}, also have other special properties. For example,
they classify the quadro-quadro Cremona transforms
(see [ESB]).

\smallpagebreak

Define the {\it tangential variety} $\tau (X)$ to be the union of all points
on all tangent stars to $X$, where the tangent star at $p_0$ is defined to
be the union of all points on all $\pp 1$'s that are the limit of secant
lines. More precisely,
${\pp 1}_*$ is in the tangent star if there exist smooth curves
$p(t),q(t)$ on $X$ such that $p(0)=q(0)=p_0$ and
$\Bbb P^1_* = \text{lim}_{t\rightarrow 0} \Bbb P^1_{pq}$.  (See [F] or
[Z1] for a more algebraic definition of the tangent star.) If we were to
require $p(t)\equiv p_0$ and only allow one point to move, we would obtain
the tangent cone.  Of course at smooth points the two notions coincide and
agree with the embedded tangent space.

For general  $X^n\subset\pp {n+a}$ one expects:
$$
\text{dim}\sigma (X) = \text{min} \{ 2n+1 , {n+a}\} \ \ \
\text{dim}\tau (X) = \text{min} \{ 2n , {n+a}\} .
$$
and $\sx$ or $\tx$ is said to be
{\it degenerate} if these fail. In this paper we will study the
local differential geometry of varieties with degenerate $\sx$ and $\tx$ and
the global geometric consequences of the local differential geometry.
Of course dim$\tx \leq$dim$\sx$,
however a much stronger statement is true:
   \proclaim {Connectnedness Applied to Secant Varieties 1.3} (Fulton-Hansen
[FL])
Let $X^n\subset\Bbb P^{n+a}$ be a variety.  Then either
$$\align &i.\ \text{dim}\,\tau (X) = 2n\ \
  \text{and dim}\,\sigma (X) = 2n+1.\\
& \ \text{or}\\
& ii.\  \sigma (X) = \tau (X).\endalign
$$
\endproclaim

The paper is organized as follows:
 In \S 2 we explain some subtleties one encounters when studying $\tx$ locally
and give some examples of varieties $X$ with   $\sx$ and $\tx$
degenerate. In \S 3 we describe the Severi varieties as the complexification
of the projective planes for the four division algebras. In \S 4 we
define the second fundamental form of a mapping and state two ways
to compute $\tdim\tx$ using
$II_{X,x}$, the second fundamental form
of $X$ at a general point $x\in X$, following
[GH].
In \S 5 we describe the second fundamental forms of the Severi varieties.
In \S 6 we describe   geometric properties of systems of quadrics that could
occur as second fundamental forms of varieties with degenerate tangential
varieties, culminating in (6.26).

In \S 7 we state a version of the rank restriction theorem
([L1],(4.14)) for varieties
with degenerate secant varieties, namely:

\proclaim{Theorem 7.1 (Rank restrictions) }
 Let  $X^n\subset\pp{n+a}$ be a smooth variety with degenerate secant
variety
 that is a hypersurface.
Let $x\in X$ be a general point.
 Let $r$ be the maximum rank of a quadric
in $|II|_x$ annhilating a $II$-generic vector in $T_xX$.
Then
$$
r\geq n-a+2.
$$
\endproclaim

Here $|II|_x$ is the system of quadrics corresponding to $II$
(see (4.12)) and a quadratic form
$Q$ is said to annhilate $v\in T_xX$ if $Q(v,w)=0$
for all $w\in T_xX$. See the notation section below for the definition
of $II$-generic.
Corollary
(7.3) provides a slight refinement of Zak's theorem on tangencies.

 In \S 8 we review some basic facts regarding moving frames
and   differential invariants of projective varieties. In \S 9
we prove the statements about how to compute $\tdim\tx$ stated in
\S 4. In \S 10 we describe how to compute $\tdim\sx$ infinitesimally.
The computation immediately implies:

\proclaim{Theorem 10.3}
Let $Y\subset \ppp V$ be a variety and let
$X=v_d(Y)\subset\ppp S^dV$ be the Veronese re-embedding.
If $d>2$ or
$d=2$ and $Y$ is not a linear subspace of $\ppp V$,
then $\sx$ is nondegenerate.
\endproclaim

\proclaim{Theorem 10.5  }  Let $X^n\subset\pp{n+a}$ be a smooth variety
not contained in a hyperplane
with degenerate secant variety of dimension $n+a_0$.
Let $\sigma_k(X)$ denote the closure of the union of
all secant $\pp{k-1}$'s of $X$. Then
$$
\text{dim}\sigma_k(X)\leq n+ (k-1)a_0.
$$
\endproclaim

(10.3) is a slight extension of a theorem of Roberts [R], where
the $d>2$ case was proved and (10.4) is
the superadditivity result of Zak-Fantechi [F].

In \S 11 we prove the statements made in \S 6. In
\S 12 we do a further local computation to prove
a refinement of ([GH], 5.17):

 \proclaim {Theorem 12.1} Let $X^n\subset\pp {n+a}$
be a
 variety (respectively a patch of a complex manifold)
 with degenerate tangential variety (resp. manifold) of
dimension
$n+a_0$. Denote the  the smooth points of
 $X$ by $X_{sm}$. Then
 the Gauss map of $\tx$ has fibers of dimension at least
$\delta_{\tau}(X)+1$, where $\delta_{\tau}(X)$
is the tangential defect of $X$.
 \endproclaim

In \S 13 we prove the key step in the rank restriction theorem (7.9) and
 we combine the rank restriction theorem and the results of
\S 12 to prove

 \proclaim {Theorem 13.10} Let $X^n\subset\pp {n+a}$
be a smooth variety with degenerate tangential variety of dimension $n+a_0$.
Then
the Gauss map of $\tx$ has
fibers of dimension at least $\delta_{\tau}(X)+2$ dimensional fibers,
where $\delta_{\tau}(X)$ is the tangential defect of $X$.
 \endproclaim

Finally, in \S 15 we do some further computations to complete
our proof of Zak's theorem on Severi varieties
(1.2).

\smallpagebreak

\noindent{\bf Acknowledgements}:
It is a pleasure to thank several people who helped the
author in the preparation of this paper:
Robert Bryant, Dennis DeTurck,  Ron Donagi, Dave Eisenbud, Joe Harris,
 Eugene Lerman,     and
and the referee. In particular,
 the referee pointed out  an error in an earlier version of (6.23).

\smallpagebreak

\noindent{\bf Notation}: By a {\it variety} we  mean the zero locus in
projective space of a collection of homogeneous polynomials
that is reduced and irreducible.  $\Bbb C\pp k$
will  be denoted $\pp k$. We will often write
$V=\bcc{k+1}$ and $\ppp V= \pp k$.
We will generally denote varieties by $X$,  the
 smooth points of $X$  will be denoted $X_{sm}$,
and the singular points $X_{sing}$.
If $Y\subset\pp m$ then $\hat Y\subset\Bbb C^{m+1}$ will denote its
deprojectivization. Alternating products   will be denoted
with a wedge ($\wedge$), and symmetric
products will not have anything (e.g. $\omega\circ\beta$ will be denoted
 $\omega\beta$).  In
general we will supress reference to the base point of our manifold $X$ when we
abbreviate the names of
fibers of bundles over $X\subset\ppp V$,
so $T$ should be read as $T_pX$ for some $p\in X$, $N$ as $N_pX$ etc... .
   $\hat T_xX\subset V$ is the
deprojectivization of the
embedded
tangent space $\tilde T_xX$.  $\{ e_i \}$ means the span of the vectors $e_i$
over the index range $i$.  If $V$ is a vector space and $W$ a subspace, and
$(e_1\hd e_n )$ a basis of $V$ such that $\{e_1\hd e_p \}=W$, we will
occasionally write
$\{e_{p+1}\hd e_n\}$ mod$W$ to denote the space $V/W$.
$\intprod $ denotes contraction, i.e. given
$P\in S^kT^*$ and $v\in T$, $(v\intprod P )\in S^{k-1}T^*$.
We will use the summation convention throughout (i.e. repeated indices are to
be summed over).   By a
$II$-{\it generic
vector} we mean a $v\in T$  based at a general point such that the
associated maps and spaces one constructs from $II$ (e.g.
$\tker II_v$ as described in \S 6)
are all locally of constant rank/dimension as $v$ is varied in $T$.

We  will use the following conventions for indicies
$$\align &1\leq A,B,C\leq n+a \\
& 1\leq \alpha ,\beta \leq n \\
& 1\leq \lambda \leq a_0 \\
& a_0+1\leq h \leq a  \\
%& 2\leq l,m \leq a_0 \\
& 2\leq \rho ,\sigma\leq n \\
& 2\leq i,j,k,l \leq r+1 \\
& r+2\leq  s,t \leq a_0 \\
& a_0+1\leq \epsilon , \delta \leq n \\
& n+1\leq \mu ,\nu \leq n+a. \endalign
 $$

 \medpagebreak

\noindent\S 2. {\bf Some remarks on secant and tangential varieties.}

  In this paper we will describe how
to compute the dimensions of the secant and tangential varieties of smooth
varieties
infinitesimally.  The tests can compute the dimensions of the secant and
tangential varieties of the smooth points of a singular variety.  However it is
not in general the case that
dim$\tau (X_{sm}) =$ dim$\tx$ (where $X_{sm}$ denotes the smooth points of
$X$).
This is because the   the tangent star at singular points is larger than the
tangent cone.
   A differential geometer might be inclined to define the tangential variety
to be the union of all tangent cones.  This would avoid having huge
tangent spaces at singular points and allow one to compute
dimensions infinitesimally
even for singular varieties.  Unfortunately, to have (1.3)
hold, we need to use the tangent star as can be seen in the
following:

\smallpagebreak

{\bf Example} (2.1):  Let $\Cal C$ be a nondegenerate smooth curve in $\Bbb
P^4$. Let
$p\in\pp 4 , p\notin \Cal C$.  Let $X= C(\Cal C,p)$ be the cone over $\Cal C$
with vertex $p$.  Then it is easy to see that $\sigma (X) = C(\sigma(\Cal C ),
p)$, so dim$\sigma (X) =$dim$\sigma (X_{sm}) =4$. The infinitesimal tests
we will develop in \S 9 and \S 10 yield
dim$\tau (X_{sm})=3$, dim$\sigma (X_{sm}) =4$. If one were to use the tangent
cone
in defining the tangential variety,
call the resulting object $\tau_{cone}(X)$,
 then one would get dim$\tau_{cone}(X)
= 3$ because except at the vertex,
$\tau (X) = C(\tau(\Cal C ), p)$, and the tangent cone at the vertex is just
$X$ itself.  But using the
 tangent star, there is a four dimensional space of tangent lines at the vertex
$p$, which allows $\tau(X)=\sigma(X)$.

\bigpagebreak

Schematic picture:

\bigpagebreak

$$ {}$$

\bigpagebreak

$$ {}\tag 2.2
$$

\bigpagebreak

$$ {}$$

For this reason we will generally restrict ourselves to studying smooth
varieties in this paper.

An example to keep in mind of a smooth variety with a degenerate secant variety
is:

\noindent{\bf Example 2.3. (The Segre $\pp{k-1}\times \pp{r-1}\ra\pp{kr-1}$)}.
 Denote the $k\times r$ matrices by $M_{kr}$.  Let $X = \Bbb P ($rank one
matrices in $M_{kr}$) and consider $X$ as a subvariety of $\pp {kr-1}$.  Here
dim$X$ = $r+k-1$ .  Then $\sigma (X) =  \Bbb P ($rank two matrices in $M_{kr}$)
and  dim $\sigma (X) = 2(r+k-1) = 2$dim$X -1$.
In some sense all known smooth varietes with degenerate secant varieties are
like these or projections  or hyperplane sections of such.

\smallpagebreak

Following [LV], for $X$ with degenerate $\sm$, the {\it secant defect} is
defined to be $\delta (X) :=2n+1 -\text{dim}\sm$.  It is difficult to find
smooth varieties with large secant  defects.  On the other hand there are
singular varieties with arbitrarily large secant defects.  One takes   $X=
\Bbb P ($rank $l$ matrices in $M_{kr}$) which has dimension $l(k+r-l)$. Then
the secant variety is $\Bbb P$(rank $2l$ matrices) which is of dimension
$2l(k+r-2l)$ and thus $\delta (X) = 2l^2 +1$.  Such $X$
 are singular
along the rank $l-1$ matrices.

\bigpagebreak

\noindent\S 3. {\bf The Severi varieties.}

For the purposes of this paper it will be better to look at the   Severi
varieties  from a unified point of view. This point of view  was first observed
independently by  J. Roberts and T. Banchoff.

Let $\Bbb A_{\Bbb R}$ denote $\Bbb R$, $\Bbb C$, $\Bbb H$, or $\Bbb O$ (where
$\Bbb O$ denotes the octonians, or Cayley numbers).  Let $\Bbb A = \Bbb A_{\Bbb
R}\ot_{\Bbb R} \Bbb C$. (From now on $\Bbb F\ot_{\Bbb R} \Bbb C$ will be
denoted by $\Bbb F^{\Bbb C}$.) Let $\Cal H_{\Bbb R}$ denote the $\Bbb A_{\Bbb
R}$-Hermitian  forms on $\Bbb A_{\Bbb R}^3$,  i.e., the $3\times 3$ $\Bbb
A_{\Bbb R}$-Hermitian matrices.   If $x\in\Cal H_{\Bbb R}$, then
we may write
$$
x = \pmatrix r_1 & \bar u_1 & \bar u_2 \\
u_1 & r_2 & \bar u_3 \\ u_2 & u_3 & r_3 \endpmatrix \ \ r_i\in \Bbb R , \
u_i\in \Bbb A_{\Bbb R}.\tag 3.1
$$
Let $\Cal H = \Cal H_{\Bbb R}\ot_{\Bbb R} \Bbb C$.

Define a cubic form $det$ on $\Cal H$ by
$$ det(x) := \frac{1}{6}((\text{trace}(x))^3 + 2\text{trace}(x^3) - 3
\text{trace}(x)\text{trace}(x^2)).
\tag 3.2
$$
The reader may check that the notion of $x^2$ and $x^3$ make sense (one needs
to use the Moufang identites (see [Hv]) in the case of the octonians).  This is
just the usual determinant of a $3\times 3$ matrix when $\Bbb A = \Bbb C$. When
$\Bbb A = \Bbb O^{\Bbb C}$ one cannot define $det$ for $4\times 4$ or larger
matrices.  The reader may find it amusing  to check what happens in the
$2\times 2$ case.

Now, considering $\Cal H$ as a vector space over $\Bbb C$, let $G$ be the
subgroup of $Gl(\Cal H,\Bbb C )$ preserving $det$, i.e., define
$$ G:= \{g\in GL(\Cal H,\Bbb C) | det(gx) = det(x)\}.
\tag 3.3
$$
The reader may check that the respective groups are:
$$\align \Bbb A_{\Bbb R} =& \ \ \ \ \ \ \ \ G = \tag 3.4\\
 &\Bbb R \ \ \ \ \ \ \ \ \ \  \ Sl(3,\Bbb R)^{\Bbb C}= Sl(3,\Bbb C)\\
  &\Bbb C \ \ \ \ \ \ \ \ \ \  \ Sl(3,\Bbb C)^{\Bbb C}= Sl(3,\Bbb C)\times
Sl(3,\Bbb C)\\
   &\Bbb H \ \ \ \ \ \ \ \ \ \  \ Sl(3,\Bbb H)^{\Bbb C}= Sl(6,\Bbb C)\\
    &\Bbb O \ \  \ \ \ \ \ \ \ \  \ \lq Sl(3,\Bbb O)^{\Bbb C}\text{\rq} = E_6
\endalign
$$
 where we have written  \lq$Sl(3,\Bbb O)^{\Bbb C}$\rq  merely to be suggestive.

$det$ tells us which elements of $\Cal H$ are of less than full rank. One can
also unambiguously define a notion of being rank one; either by taking $2\times
2$ minors or by noting  that under the $G$ action each $x\in \Cal H$ is
diagonalizable and one can take as the rank of $x$ the number of nonzero
elements in the diagonalization of $x$.

Let
$$
X := \pp {}\{\text{ rank one elements of }\Cal H \} = \pp {}\{ G\text{ orbit of
any rank one matrix}\}.
\tag 3.5
$$
It is not difficult to check that $X = (\Bbb A_{\Bbb R}\pp 2)^{\Bbb C}$.
The four Severi varieties  are  exactly $X\subset \pp{}(\Cal H )$  for
$\Bbb A_{\Bbb R} = \Bbb R,\Bbb C, \Bbb H,\Bbb O$ (where again we treat
$\Cal H$ as a  vector space over $\Bbb C$).  Note that
$\sm = \pp{}\{\text{rank two elements}\}\neq\pp {} (\Cal H )$.

\bigpagebreak

\noindent\S 4 {\bf Second Fundamental Forms}.

It is possible to define second fundamental forms in quite general
settings. For the purposes of this paper, it will be useful to
deal with the second fundamental form of a mapping to a projective space.

Let $Y$ be a   (smooth) complex manifold and let
$\phi : Y\ra \ppp W$ be a smooth holomorphic mapping.
(e.g. $\phi$ could be the restriction of a rational mapping to a neighborhood
of a smooth regular point). Fixing $y\in Y$, one has the
{\it differential } of $\phi$ at $y$;
$$
\phi_{*y} : T_yY\ra T_{\phi (y)}\phi (Y)\subset T_{\phi (y)}\ppp W\tag 4.1
$$
and the {\it Gauss map} of $\phi$, (where   $\tdim\phi(Y)=k$);
$$
\align \gamma_{\phi}: &Y\ra G(k+1,W)\tag 4.2\\
& y\mapsto \hat T_{\phi (y)}\phi (Y) =: \hat T\endalign
$$
where $\hat T$ is the cone over the embedded tangent space
and $G(k+1,W)$ denotes
the Grassmanian of $(k+1)$-planes in $W$.
Taking the derivative of $\gamma_{\phi}$, at $y$, one
obtains a linear map
$$
\gamma_{\phi * y} : T_yY\ra T_{\hat T} G(k+1,W) = \hat T^* \ot W/\hat T
\tag 4.3
$$
Let $\hat\phi(y)\subset \hat T$ denote the line corresponding to $\phi (y)$.
Observe that
$$
\align &   \forall v\in T_yY,\ \hat\phi (y) \subset \tker \gamma_{\phi *y}(v)
\tag 4.4\\
&  \tker \phi_{*y}\subseteq \tker \gamma_{\phi * y}.\tag 4.5\endalign
$$
Thus $\gamma_{\phi *y}$ descends to an element of
$(T_yY/\tker\phi_{*y})^*\ot (\hat T/ \hat\phi(y))^*\ot W/\hat T$.
We may identify $T_yY/\tker\phi_{*y}\simeq \phi_{*y}(T_yY)\simeq
T_{\phi(y)}\phi(Y)\subseteq T_{\phi(y)}  \ppp W$.
Moreover, since
$$
\align &
T^*_{\phi(y)}\phi(Y)=(\hat T/{\hat\phi(y)} )^*\ot \hat\phi(y)\tag 4.6\\
&N_{\phi(y)}\phi(Y)=\hat\phi(y)^*\ot W/\hat T\endalign
$$
we obtain an element of
$T^*_{\phi(y)}\phi(Y)\ot T^*_{\phi(y)}\phi(Y)\ot N_{\phi(y)}\phi(Y)
$
and, in fact, because the Gauss map is already the derivative of
a map and mixed partials commute, we obtain
$$
II_{\phi , y}\in S^2(T_{y} Y/\tker\phi_{*y})^*\ot N_{\phi(y)}\phi(Y)\tag 4.7
$$
the {\it second fundamental form} of the mapping $\phi$ at $y$.
If $\phi$ is the inclusion map, we recover the second fundamental
form of a subvariety.

\smallpagebreak

Let $X^n\subset\pp{n+a}$ be a variety and let $x\in X$ be a general point.
Let $II=II_{X,x}$ denote the second fundamental form at $x$. We will
say a vector $v\in T=T_xX$ is {\it $II$-generic} if all the associated
maps and spaces one constructs from $II$ and $v$ (see e.g. $II_v$
below) are locally of constant rank/dimension near $v$.

Fix a $II$-generic $v\in T$ and consider the map:
$$
\align II_v : T &\ra N \tag 4.8\\
 w&\mapsto II(v,w)\endalign
$$
and the rational map:
$$
\align ii : \ppp T &- - \ra \ppp N \tag 4.9\\
 [w]\ \ &\mapsto [II(w,w) ].\endalign
$$

\proclaim{Proposition 4.10}Let $X^n\subset\pp{n+a}$
be a smooth variety, let $x\in X$ be a general point and let
$v\in T_xX$ be a $II$-generic vector. Then
$$
\align
\tdim\tx &= n + \tdim II_v(T) \tag 4.11\\
&= n + \tdim ii(  \ppp T).\endalign
$$
\endproclaim
We give a proof of the first line of (4.11) in \S 9
(originally due to Terracini, [T] and rephrased in modern language
in [GH]) and prove
  $\tdim II_v(T)=\tdim ii(  \ppp T)$ in (6.6).

Considering the second fundamental form as a mapping
$II^* :N^*\ra S^2T^*$, let
$$
\ii = \ppp ( II^*(N^*)) \subset \ppp (S^2T^*).
\tag 4.12
$$

\bigpagebreak

\noindent \S 5. {\bf Second Fundamental Forms of the Severi Varieties}

 Let $Y\subset \pp {n-1}$ respectively denote $ \emptyset ,\ \pp 1 *\pp 1
, \ \pp 1\times\pp 3  , \ \Bbb S$,  where respectively $n= 2,4,8,16$;  $*$
denotes disjoint union, $\pp 1\times\pp 3$ occurs as the Segre, and $\Bbb S $
is the 10 dimensional spinor variety obtained by taking the projectivized orbit
of a highest weight vector in $\Bbb C^{16}$ under the $Spin(10,\Bbb C)$ action.

 The four Severi varieties have:
$$
|II|=\text{ the complete system of quadrics vanishing on }Y.\tag 5.1
$$
 If we write a
Severi variety as $G/P$, the fastest way to see this is  to observe that
  $Y$ must be an orbit of the $P$ action on $T_x(G/P)$.  Now dim$|II| =
$codim$(G/P)$ since $II$ must surject onto $N$ (by (10.2)), so $Y$ must be
small enough to have an $\frac{n}{2} + 2 $ dimensional space of quadrics
vanishing on  it.  But in each case here (except $n=2$ which is easy) $Y$
is the projectivized orbit of a highest weight vector and therefore the
smallest possible (nonempty) orbit.  Alternatively, one may compute directly,
see [GH] for the classical cases and the last chapter of [Hv] for  a good model
to study the $E_6/P$ case  (actually [Hv] only deals with $Spin(9,\Bbb R)$ on
$\Bbb R^{16}$ but the argument for the case here is similar).

\smallpagebreak

  Another way to describe the second fundamental form is to use local affine
coordinates.  If we want to know the second fundamental form at $p$, we choose
coordinates $[X_0\hd X_N]$ of $\pp N$ such that $p = [1,0\hd 0]$, and
 affine coordinates
$x^{B}= X_B/X_0$ where $X$ is locally describable as a graph over
$x^1\hd x^n$  and $T_p \simeq \{x^{\alpha}\}$,
$1\leq \alpha ,\beta\leq n$.  In this case in some
neighborhood
of $(0\hd 0)$,
$$
x^{\mu} = x^{\mu}(x^{\alpha})= \qq {\mu}{\alpha}{\beta}x^{\alpha}x^{\beta}
+\text{terms of order three or greater}
$$
where
$n+1\leq\mu\leq \na$, and
 $\qq {\mu}{\alpha}{\beta}$ are the coefficients of the second fundamental
form (see [GH] for more details).

 In the case of the Severi varieties, choose affine coordinates based at $[p]$
where
$$
p = \pmatrix 1& 0 & 0\\
 0 & 0& 0\\
   0 & 0& 0\endpmatrix
$$
and following the notation of (3.1), denote the affine coordinates
$u_1,u_2,u_3\in\Bbb A ,r_2,r_3\in \Bbb C$ where the tangent space to $p$ is
$\{u_1,u_2\}$ (the span is taken over $\Bbb C$).  Here we have (where $\equiv$
means modulo terms of order three or greater)
$$
r_2(u_1,u_2)\equiv u_1\bar {u_1}\ \text{as} \  \ \ \text{det}\pmatrix
1&\bar {u_1}\\ u_1 & r_2\endpmatrix =0
$$
$$r_3(u_1,u_2)\equiv  u_2\bar {u_2}\ \text{as} \  \  \ \text{det}\pmatrix
1&\bar {u_2}\\ u_2 & r_3\endpmatrix =0
$$
$$u_3(u_1,u_2)\equiv \bar {u_2}u_1\ \text{as} \  \  \ \text{det}\pmatrix
1&\bar {u_1}\\ u_2 & u_3\endpmatrix =0
$$
where the last equation gives us one, two, four or eight quadratic forms.   The
determinants come from the vanishing of $2\times 2$ minors that must be zero to
make the Severi variety consist only of rank one elements.  (In fact,
if one checks the other $2\times 2$ minors, one sees that the equations above
are exact, i.e. all third and higher order terms vanish. The idea of the proof
of Zak's theorem on Severi varieties given in this paper is to first show that
any putative Severi variety must have
one of the four fundamental forms in (5.1) and then to show  all third  and
higher
order invariants   must be zero.)

In division algebra notation the second fundamental forms
are
$$
  |II| =  \ppp \{ u_1\bar{u_1}, u_2\bar {u_2},\bar {u_2}u_1\}. \tag 5.2
$$

 For example if $\Bbb A
=
\Bbb O^{\Bbb C}$ and we write
$u = (u^0\hd u^7 )= u^0 + u^1J_1 +\hdots + u^7J_7$ then
$$
 \align
 |II| = &
\ppp \{ (u^0_1)^2 +\hdots + (u^7_1)^2,\ (u^0_2)^2 +\hdots + (u^7_2)^2,\
u^0_1u^0_2+\hdots u^7_1u^7_2, \\
& \ \  u^0_1u^1_2-u^1_1u^0_2 + u^2_1u^3_2-u^3_1u^2_2 +\hdots + u^7_1u^4_2,
\hdots , u^0_1u^7_2 +\hdots - u^4_1u^1_2 \}\endalign
$$
where we have used the octonionic multiplication table:

\bigpagebreak

$$ {}$$

\bigpagebreak

\bigpagebreak

$$ {}$$

\bigpagebreak

$$ {}$$
where if $J_1\hd J_7$are a basis of Imag($\Bbb O^{\Bbb C}$) then
the multiplication is to be read cyclicly and counter clockwise, e.g. $J_1J_2 =
-J_2J_1=J_3,J_2J_3 = J_1,J_1J_7 = J_4,J_4J_2 = -J_6$.

\smallpagebreak

   In fact, the original description of the Severi varieties used by Zak to
prove his theorem was to consider $Y\subset \pp{n-1}$,   and linearly embed
   $$
   Y\subset \pp{n-1}\subset\pp n
   $$
and to note that the image of the natural embedding by the linear system of
quadrics on $\pp n$ vanishing along $Y$ of $Bl_{Y}\pp
n\rightarrow
\pp{\frac{3n}{2} +2}$
is   the Severi variety.

\bigpagebreak

\noindent\S 6. {\bf Properties of systems of quadrics with a
(tangential) defect}

Let $T$ be an $n$-dimensional vector space,
$N$ an $a$ dimensional vector space with $a\leq n$, and let
$II\in S^2T^*\ot N$.
We will write $II :S^2T\ra N$ or
$II^* :N^*\ra S^2T^*$ when we want to consider
$II$ as a linear map.
Assume $II^*$ is injective, which will be the case in
our situation by (10.2).  Let $\ii = \ppp
(II^*(N^*))\subseteq
\ppp S^2T^*$, and let $\hii =II^*(N^*)\subseteq   S^2T^*$.
  Given $Q\in \hii$, let $[Q]$ denote the
hypersurface in $\ppp T$ it determines. If  $v\in T$ and
$v\intprod Q=0$, i.e. $[v]\in [Q]_{sing}$, we will say
  $Q$ {\it annihilates} $v$.

We will say $II$ has a {\it (tangential) defect}
if for every $v\in T$ there is a $Q\in \hii$ annihilating $v$.
If $X^n\subset\pp{n+a}$ is a variety with degenerate tangential
variety, $x\in X$ is a general point and $a\leq n$, then
$II_{X,x}$ has a defect (4.10). In what follows, we will fix a
$II$-generic $v\in T$, which is natural, since in the geometric
setting this
amounts choosing a smooth point of $\tx$ in $\tilde T_xX$
up to scalings.

\smallpagebreak

Fixing a $II$-generic $v\in T$, consider the following
linear subspace of $\ii$:
$$
\align
\tann (v) &= \{ [Q]\in\ii\ |\ [v]\in [Q]_{sing} \}\tag 6.1\\
&= \ppp \{ Q\in \hii\ |\ v\intprod Q=0 \}\endalign
$$
the annhilator of $v$. Also consider
$$
\align
\tsingloc (\tann (v)) &=
\{ v\in T \, | \, [v]\in [Q]_{sing}\ \forall\ [Q]\in\tann (v) \}\tag 6.2 \\
&= \{ v\in T \, | \, v\intprod Q=0\
\forall \ Q\in\widehat\tann (v) \},
\endalign
$$
the singular locus of the zero set of $\tann (v)$.
Consider the linear map
$$
\align II_v :T&\ra N\tag 6.3\\
w&\mapsto II(v,w).\endalign
$$
Note   that
$$
\ppp (II^* (II_v(T)\upperp ))=\tann (v).\tag 6.4
$$

\proclaim{Lemma 6.5}
   $\ppp \{ v,\tker
II_v \}\subseteq  \tsingloc (\tann (v))$
\endproclaim

 We will give two proofs of (6.5), one in frames in \S 11 and
one that will follow from the following:

\proclaim{Proposition 6.6} Consider the rational map
$$
\align ii :\ppp T -& - \rightarrow \ppp N\\
 [v] &\mapsto [II(v,v)]. \endalign
$$
Let $v\in T$ be $II$-generic.

1. Let  $w\in T$, and let
$\underline w= v^*\ot (w  \tmod  v)$, so
$\underline w\in T_{[v]}\ppp T$,
then
$$
ii_{*[v]}(\underline w)= II(v,v)^*\ot (II(v,w)\tmod \{ II(v,v)\} ).\tag 6.7
$$
In particular,
$$
\tker (ii_{*,[v]})= v^*\ot (\tker II_v \
\tmod\ \{ v \} )\tag 6.8
$$

2. For all $w_1,w_2\in T$, let
$\tilde w_j =
v^*\ot (w_j\tmod \{ v,\tker II_v \})$, so $\tilde w_j\in
T_{[v]}\ppp T/
\tker ii_{*[v]}$, then
$$
II_{ii,[v]}(\tilde w_1 ,\tilde w_2 )
= II(v,v)^*\ot ( II(w_1,w_2 ) \ \tmod \ II_v(T)). \tag 6.9
$$
where  $II_{ii,[v]}$ is the second fundamental form of the mapping
$ii$
at $[v]$ as in (4.7).

In particular,  $\tsingloc |II_{ii,[v]}|=
v^*\ot (\tsingloc (\tann (v)) \tmod  \{ v,\tker
(II_v) \} )$
\endproclaim

Note that (6.8) proves the equivalence of the two lines in (4.11).

\demo{Proof of (6.6)}
Let $v_t$ be a curve with $v_0=v$. We have
$$
ii(v_t) = II(v_t,v_t) \tag 6.10
$$
so
$$
 \frac d{dt}ii(v_t)|_{t=0} \equiv 2II(v,v_0') \ \tmod\ II(v,v) \tag 6.11
$$
where $v_0'=\frac d{dt} v_t|_{t=0}$, proving 1. To see 2,
we compute $II_{ii,[v]}$ as in [GH];
$$
\align
 \frac {d^2}{(dt)^2}ii(v_t) |_{t=0} &\equiv 2II(v,v_0'') + 2II(v_0',v_0')
\ \tmod II_v(T) \\
&\equiv
2II(v_0',v_0')
\ \tmod II_v(T) \tag 6.12 \endalign
$$
polarizing (6.12) and putting in twists,
one recovers $II_{ii,[v]}$.
\qed\enddemo

\demo{Proof of (6.5)} Apply (4.5) with $\phi = ii$
and use (6.8) and (6.4).\qed
\enddemo

(6.5) is a statement about a
  special property
of the subsystem of
 quadrics in $\ii$  annhilating  a fixed vector.
In general, one cannot expect to make any statements about
all the quadrics in $\ii$, as it may be the case that a subsystem
of $\ii$ already has a defect
of the same size and therefore only this subsystem will have
special properties.
What we really want to study is the smallest subsystem
of $\ii$
generated by quadrics annhilating the
$II$-generic vectors.

Let $Z$ denote the (closure of the) image of the rational map $ii$.
It may be the case that $Z$ is a cone. If so, its vertex is
$$
V_Z = \cap_{z\in Z_{sm}} \tilde T_zZ \subset Z. \tag 6.13
$$
($Z_{sm}$ denotes the
smooth points of $Z$.) Note that we may take the intersection over
 $\{ z\in Z_{sm}
| z=ii([v])\text{ for some }v\in T^{gen} \}$.
Here and in what follows,
 $T^{gen}$ denotes the space of $II$-generic vectors.

\proclaim{Proposition 6.14}   $II^* (\hat V_Z\upperp )$ is the smallest
subsystem of $\hii$ with the same size
tangential  defect.
\endproclaim

\demo{Proof}
By (6.7), $\hat T_{ii([v])}Z =   II_v(T)\subset N$.
So,
$\widehat\tann (v) = II^*(II_v(T)\upperp)=
II^* ((\hat T_{ii([v])}Z)\upperp )$
and therefore
$$
\text{ span}_{v\in T^{gen}}\widehat\tann (v) = II^*(
 (\cap_{v\in T^{gen}}\hat
T_{ii([v])}Z) \upperp ).\qed
$$
\enddemo

 The remainder
of the results in this section can be proved without using moving
frames, but already our proof of (6.16) without frames was quite messy
so we will only present proofs in frames
and defer them to \S 11. Readers particularly
allergic to frames may wish to write their own proofs of what follows.
See [D] for a frame free  argument of a slightly incorrect version of (6.16).
(One mistake is that he   implicitly assumes $\gamma_Z$
is nondegenerate.)

Let
$$
F_v=\gamma_Z\inv (\gamma (ii([v])))\tag 6.15
$$
i.e., $F_v$ is the fiber of the Gauss map of $Z$ through
$ii([v])$.

\proclaim{Lemma 6.16}  $\tsingloc (\tann (v) )\subseteq
\tbaseloc  II^*(F_v\upperp )\ .$
\endproclaim

\smallpagebreak

By (6.5) we may consider the isomorphism
$$
\underline{II}_v :
T/\tsingloc (\tann (v)) \ra II_v(T)/
 \hat F_v. \tag 6.17
$$
By (6.16), for all $w\in \tsingloc (\tann (v))$ we have
$$
II_w (\tsingloc (\tann (v)))
\tmod \hat F_v
 \subseteq II_v(\tsingloc (\tann (v)))
\tmod \hat F_v \tag 6.18
$$
so the maps
$$
II_w'
: T/\tsingloc (\tann (v))\ra II_v(T)/
 \hat F_v   \tag 6.19
$$
are well defined. In fact it will be useful to consider these as
endomorphisms of
$T/\tsingloc (\tann (v))$. To this end, for $w\in \tsingloc (\tann (v))$, let
$$
\phi_w  = (\underline {II}_v)\inv\circ II_w'
 : T/\tsingloc (\tann (v))\ra T/\tsingloc (\tann (v)) \tag 6.20
$$

If $Z$ is a hypersurface, then
 $T/\tsingloc (\tann (v))$ comes equipped with a
quadratic form defined up to scale, namely $\tann (v)$. We will show that
in this case, for
$w\in \tker II_v$ that
 $\phi_w\in\frak c\frak o (T/\tsingloc (\tann (v)),\tann
(v))$. (Note that by construction, $\phi_v = Id_{T/\tsingloc (\tann (v))}$.) In
fact it will be useful to choose a quadratic form
$P_v\in \widehat{\tann (v)}$ so we do so.

\proclaim{Lemma 6.21} If $Z$ is a hypersurface,
then for all $w\in \tker II_v$,
 $$
\phi_w\in\frak s\frak o (T/\tsingloc (\tann (v)), P_v).
$$
\endproclaim

In general, if a variety $Y$ is a cone with vertex $V_Y$, then the fibers of
$\gamma_Y$, the Gauss map of $Y$, all contain $V_Y$. If $Y$ is a cone over a
variety
having a nondegenerate Gauss map, then the dimension of the
general fiber of $\gamma_Y$ has dimension equal to $\tdim V_Y+1$.
This will be the case if  and only if $Y$ is not
\lq\lq built out of a tangent developable\rq\rq , in the language
of ([GH], 2.27).
Recall that  a variety has a degenerate Gauss map
if and only if it is
 ruled by linear spaces along which the embedded tangent
space is constant. Two types of varieties having this
property are cones and tangent developables. If
one takes a cone over a tangent developable or visa versa,
the resulting variety will still have a degenerate Gauss map,
and is an example
of a variety \lq\lq built out of a tangent developable\rq\rq .

{}From now on, assume $Z$ is a hypersurface and that
$\tdim V_Z+1  = \tdim F_v$, i.e., that
$Z$ is not built out of a tangent developable.
By (6.16),
$$
 II^*(V_Z\upperp)|_{ \{ v,\tker II_v \}} \tag 6.22
$$
is a well defined  quadric up to scale. Let
 $Q_v\in II^*({V_Z}\upperp)|_{ \{ v,\tker II_v \}}$
 be the unique quadric such that $Q_v(v,v)=1$.

\proclaim{Lemma 6.23}
Assume $Z$ is a hypersurface
such that  $\gamma_Z$
has fibers of dimension $\tdim V_Z-1$.
(This condition is automatic if $\gamma_Z$ is nondegenerate.)
For all $w_1,w_2\in \{v , \tker II_v \}$,
$$
\phi_{w_1}\circ\phi_{w_2}
+ \phi_{w_2}\circ\phi_{w_1} + 2Q_v(w_1,w_2) Id_{T/\tsingloc (\tann (v))}=0
\tag 6.24
$$
where
$\phi_w$ is as in (6.20),
$Q_v$ is the quadric defined above,
and  $Id_{T/\tsingloc (\tann (v))}$ is the identity map on
 $T/\tsingloc (\tann (v))$. In particular,
$T/ \tsingloc (\tann (v))$ is a $Cl( \{v, \tker II_v \}, Q_v)$ module, where
$Cl( \{v, \tker II_v \}, Q_v)$ denotes the Clifford algebra.
\endproclaim

\demo{Proof of \lq \lq In particular\rq\rq}
The fundamental lemma of Clifford Algebras ([Hv] or [LM]) says:  Given a linear
map $\phi : V\rightarrow \bold A$ from a vector space $V$ with
 inner product
$Q$ into an associative algebra $\bold A$ with unit, if
$$
\phi (x)\phi (y) +\phi (y)\phi (x) + 2 Id Q(x,y) = 0\tag 6.25
$$
 then $\phi$ has a unique extension $\tilde\phi : Cl(V, Q) \rightarrow \bold
A$.
 Now let $V = \{ v,\tker II_v \}$, $Q=Q_v$,
 and $\bold A =$End$(T/\tsingloc (\tann
(v)))$.
 \qed
\enddemo

In the geometric setting, (6.23) implies

\proclaim{Theorem 6.26 (The Clifford algebra structure)}
Let $X^n\subset\Bbb C\pp{n+a}$ be a   variety with degenerate
tangential variety that is a hypersurface. Let $x\in X$ be
a general point and let $v\in T_xX$ be a $II$-generic vector.
Assume that the hypersurface
$Z\subset\ppp N_xX$  defined by $II_x$ as above (6.13)  is
such that  the Gauss map of $Z$, $\gamma_Z$,
has fibers of dimension $\tdim V_Z-1$,
where $V_Z$ denotes the vertex of $Z$ as in (6.13).
(This condition is automatic if $\gamma_Z$ is nondegenerate.)
Then the space $\{ v,\tker II_v \} $ comes equipped with a
quadratic form $Q$  and
$T_xX/ \tsingloc (\tann (v))$ is a
 $Cl(\{ v,\tker II_v \},Q)$
module, where
$Cl(\{ v,\tker II_v \},Q)$ denotes the Clifford algebra.
\endproclaim

In the critical codimension case of $a=\frac n2+2$, the vertex dimension
hypothesis in (6.26) is  satisfied when $X$ is smooth. One might
hope to show the vertex
dimension hypothesis in (6.26) always must hold
when $X$ is smooth using (7.1) below. Such a condition would have
consequences for the admissible values of tangential defects.

\bigpagebreak

\noindent\S 7. {\bf Rank Restrictions}

\proclaim{Theorem 7.1 (Rank restrictions) }
 Let  $X^n\subset\pp{n+a}$ be a smooth variety with degenerate secant
variety
 that is a hypersurface.
Let $x\in X$ be a general point.
 Let $r$ be the maximum rank of a quadric
in $|II|_x$ annhilating a $II$-generic vector.
Then
$$
r\geq n-a+2.
$$
\endproclaim

Since $\tdim (\tsingloc (\tann (v))) =n-r$, and $\tdim (\tker II_v)=
n-a+1$, it follows that
$$
\align
\tdim {F_v}
 &= \tdim (\tsingloc (\tann (v))) /\tker II_v  -1\tag 7.2\\
&= a-r-2\endalign
$$
where $F_v$ is as in (6.15)

\proclaim{Corollary 7.3} $X^n\subset\pp{n+a}$ be
 a smooth variety with degenerate
secant variety that is a hypersurface.
Let $x\in X$ be a general point and $v\in T_xX$ be a $II$-generic
vector, and  let $F_v$ be as in (6.15).
 Then
$$
a \geq \frac n2 + 2 + \frac 12 \tdim {F_v}. \tag 7.4
$$
\endproclaim

(7.3)   gives  a slightly refined version of  (1.1)
because one can always project to the case $\sx$ is a hypersurface.

\demo{Proof of (7.1)}   Consider the incidence correspondences:
$$
\align &\ci := \{ (x,H) | \tilde T_xX\subset H\} \subset X\times X^*\tag
7.5\\
& \\
& \\
X \ \ \ \ \ & \ \  \ \ \  X^*\endalign
$$
$$
\align &\Delta := \{ (x,H) |
\tilde T_xX\subset H \text{ and } Q_H=\tann (v)
\text{ for some } v\in T^{gen}_xX\}\subseteq
  \ci\tag 7.6\\
& \\
& \\
X \ \ \ \ \ & \ \  \ \ \  X^*_{\Delta}\endalign
$$
where $X^*$ is the dual variety of $X$,
 $Q_H\in  \ii_x$   corresponds to
the
image under $II^*$ of the line in $N^*$ determined by $H$,
$T^{gen}$ is the space of $II$-generic vectors, and
$X^*_{\Delta}\subseteq X^*$ is defined to be $\rho_{\Delta}(\Delta )$.
   Since $\sm$ is degenerate, we can apply  Zak's
theorem on tangencies [Z1]
 to the image of a projection of $X$ to find that
$$
a-2\geq \text{dim }\pi (\rho \inv (H))
\geq
\text{dim }\pi_{\Delta}(\rho_{\Delta}\inv (H)). \tag 7.7
$$
 Take $H$ such that $(x,H)\in \Delta$
is a general point of $\Delta$ in the sense
that both $\pi_{\Delta}$ and $\rho_{\Delta}$ are smooth at
$(x,H)$ and of maximal rank, and $x$ is a general point of $X$.

Note that
$$
\tdim\Delta = n+r,\tag 7.8
$$
where by this we mean the local dimension near $(x,H)$,
because\linebreak
 $\tdim\ppp ( II^*(\hat F_v\upperp ))=r$.
($II^*$ is injective by (10.2).)

We claim that
$$
\tdim X^*_{\Delta}=2r \tag 7.9
$$
(where again, we mean the local dimension near $H$).
To understand (7.9) heuristically, recall that $\tdim X^* = r'+(a-1)$
where  $r'$ is the rank of a generic quadric in the second fundamental
form of $X$ at a general point. The way to think of $\tdim X^*$
 is
that there are $a-1$ \lq\lq vertical\rq\rq\ contributions to the dimension
(all the hyperplanes tangent to a point) and
$r'$ \lq\lq horizontal\rq\rq\
contributions (where if $[Q_H]\in\ii$ is the quadric corresponding
to a tangent  hyperplane $H$, then $r'=\tdim (T/\tsingloc (Q_H))$). The same is
true  for $X^*_{\Delta}$, there
is a contribution of  $r$ \lq\lq vertical\rq\rq\
 directions
(corresponding to $\ppp II^*(\hat F_v\upperp )$ )
 and $r$ \lq\lq
horizontal\rq\rq\ ones
corresponding to $T/\tsingloc (\tann (v))$.
   (7.9) will be proved   in \S 13 using moving frames.

Putting everything together, we have
$$
a-2
 \geq \tdim  \pi_{\Delta}(\rho_{\Delta}\inv (H)) = \tdim \Delta - \tdim
X^*_{\Delta}
 = (n+r) -2r = n-r.\qed \tag 7.10
$$
\enddemo

\noindent\S 8. {\bf Moving frames}

We will briefly review some notation here, the reader  is refered to [L1], [L2]
for more details.
Let $X^n$ be a (patch of  a) smooth submanifold of $\ppp V=\pp{n+a}$
 and let
$x\in X$.
 Let $\cf^1$ be the bundle of all frames
$f=(\aa 0\hd\aa{n+a})$
(i.e. bases of $V$) adapted to the filtration $\hat
x\subset \hat T\subset V$,
  where $\hat x\subset V$ is the line determine by $x$ and
$\hat T $
is the deprojectivization of the embedded tangent space to $X$ at $x$.
 Write  $f=(\aa 0 ,\aa\alpha ,\aa\mu )$ for an element
of $\cf^1$, where $1\leq\alpha , \beta\leq n,\  n+1\leq\mu , \nu\leq n+a$
where $[\aa 0 ]=x$, and  $\{\aa 0\hd\aa n\}=\hat T$.
The motions in the fiber over  $[\aa 0 ]$ are $f\mapsto fg$ where
$$
g = \pmatrix
\gg 0 0 & \gg 0 \beta & \gg 0\nu \\
 0 &\gg\alpha\beta & \gg\alpha\nu \\
 0 & 0 & \gg\mu\nu\endpmatrix , \ \ \text{det}g\neq 0.\tag 8.1
$$
Let
$\Omega = f\inv df$ denote the Maurer-Cartan form.  We   write $\Omega$
as
 $$
\Omega = \pmatrix \ooo 0 0 & \ooo 0 \beta & \ooo 0 \nu \\
\oo\alpha  & \ooo \alpha\beta & \ooo \alpha\nu \\
0 & \ooo\mu\beta & \ooo\mu\nu \endpmatrix \tag 8.2
$$
where the entries of
$\Omega$ are one-forms on $\cf^1$ and we use the notation $\oo B$ for $\ooo B
0$. The Maurer-Cartan equation $$d\Omega = -\Omega\ww\Omega\tag 8.3
$$
will be used repeatedly in what follows. For example,  using (8.3) and
 the Cartan lemma applied
to $0 = d\oo\mu$
implies
$\ooo\mu\beta = \qq\mu\beta\alpha\oo\alpha$ with $\qq\mu\alpha\beta =
\qq\mu\beta\alpha$. This   gives an alternative definition of
the  second fundamental
form because
$$
  \ooo\mu\beta\oo\beta \ot A_{\mu}\text{mod }\hat T
 =\qq\mu\alpha\beta\oo\alpha\oo\beta \ot A_{\mu}\text{mod }\hat T\tag 8.4
$$
descends  (after twisting by $\aa 0^*$)
to a well defined section of $ S^2T^*\ot N$ over $X$ which is $II$.

Similarly, the quantity
$$
 \ooo\nu\mu\ooo\mu\beta\oo\beta \ot A_{\nu}\text{mod }\hat T\up 2 \tag 8.5
$$
 (where $\hat T\up 2= \{ \hat T , \widehat{II(S^2T)}\}$)
descends after twisting to  a well
defined section of $ S^3T^*\ot (N/ II(S^2T))$ over $X$,
called the {\it third fundamental form} and denoted $III$.

Following [GH], for varieties with degenerate tangential varieties,
or varieties where $a>n$,
one can define a refinement of $III$. Namely, fixing a
$II$-generic vector $v\in T$, define
the {\it refined third fundamental form},
$$
  III^v  \in S^3 (\tsingloc (\tann (v)))^*\ot N/II_v(T)\tag 8.6
$$
  by noting the quantity
$$
\ooo\nu\mu\ooo\mu\beta\oo\beta|_{S^3(\tsingloc (\tann (v)))}
\ot A_{\nu}\text{mod
}\{\hat T, \widehat{II_v(T)}\}\tag 8.7
$$
descends (after twisting) to be well defined over $X$.
   $III^v$ is well defined as if $w\in
\tsingloc (\tann (v))$, then $II_w(T)\subseteq II_v(T)$.
  One can
recover $III$ from knowing $III^v$ for all $v\in T$.

Differentiation (via (8.3))
of the equation
$$
\ooo\mu\beta-\qq\mu\beta\alpha\oo\alpha = 0\tag 8.8
$$
provides another third order invariant,
which has been called the  {\it cubic
form}.
$$
\align
 \pii & = (-d\qq\mu\alpha\gamma-\qq\mu\alpha\gamma\ooo 00 +
\qq\mu\beta\gamma\ooo\beta\alpha + \qq\mu\alpha\beta\ooo\beta\gamma -
\qq\nu\alpha\gamma\ooo\mu\nu )\oo\alpha\oo\gamma \ot (A_{\mu}
\text{ mod}\hat T)\tag 8.9\\
&=\rr\mu\alpha\beta\gamma\oo\alpha\oo\beta\oo\gamma
\ot (A_{\mu}
\text{ mod}\hat T).\endalign
$$
The Cartan Lemma  is what allows us to write the second line, i.e.
that there are functions $\rr\mu\alpha\beta\gamma$, symmetric in
their lower indices, such that
$$
\rr\mu\alpha\beta\gamma\oo\gamma=
-d\qq\mu\alpha\beta-\qq\mu\alpha\beta\ooo 00 -
\qq\mu\beta\gamma\ooo\gamma\alpha +  \qq\mu\alpha\gamma\ooo\gamma\beta -
\qq\nu\alpha\gamma\ooo\mu\nu.
\tag 8.10
$$
Although $\pii$ is well defined on $\cf^1$ it varies on the fiber as follows:
If
we move along the
fiber by motions
$$
A_{\beta}\ra A_{\beta} +g^0_{\beta}A_0 \tag  8.11
$$
 then
 the expression (8.9) will change
by
$$
 \pii\ra \pii + \frak S_{\alpha\beta\gamma}
g^0_{\beta}\qq\mu\alpha\gamma
\oo\alpha\oo\beta\oo\gamma\ot (A_{\mu}
\text{ mod}\hat T) \tag 8.12
$$
(where $\frak S$ denotes cyclic summation).  If we move in the fiber by motions
$$
A_{\mu}\ra A_{\mu} +g^{\delta}_{\mu}A_{\delta} \tag 8.13
$$
then
 the expression (8.9) will change
by
$$
 \pii\ra \pii - \frak S_{\alpha\beta\gamma}
g^{\delta}_{\nu}\qq\mu\alpha\delta\qq\nu\beta\gamma
\oo\alpha\oo\beta\oo\gamma\ot  (A_{\mu}
\text{ mod}\hat T). \tag 8.14
$$

If we differentiate    (8.10), we obtain a fourth order invariant,
$$
\align
\partial^2II &=
(-d\rr\mu\alpha\beta\gamma -2\rr\mu\alpha\beta\gamma\ooo 0 0
-\rr\nu\alpha\beta\gamma\ooo\mu\nu + \frak
S_{\alpha\beta\gamma}\rr\mu\alpha\beta\delta\ooo\delta\gamma \\
& \ \ \
-\frak S_{\alpha\beta\gamma} \qq\mu\alpha\delta\qq\nu\beta\gamma\ooo\delta\nu -
\frak S_{\alpha\beta\gamma}\qq\mu\alpha\beta\ooo 0\gamma
)\oo\alpha\oo\beta\oo\gamma
\ot\aa\mu\tmod \hat T\tag 8.15\\
&=\rr\mu\alpha\beta{\gamma\delta}\oo\alpha\oo\beta\oo\gamma\oo\delta
\ot\aa\mu\tmod \hat T\endalign
$$
where again the Cartan lemma says that there are functions
$ \rr\mu\alpha\beta{\gamma\delta}$, symmetric in their lower
indices, having the property that
$$
\align
 \rr\mu\alpha\beta{\gamma\delta}\oo\delta &=
-d\rr\mu\alpha\beta\gamma -2\rr\mu\alpha\beta\gamma\ooo 0 0
-\rr\nu\alpha\beta\gamma\ooo\mu\nu + \frak
S_{\alpha\beta\gamma}\rr\mu\alpha\beta\delta\ooo\delta\gamma\tag 8.16 \\
& \ \ \
-\frak S_{\alpha\beta\gamma} \qq\mu\alpha\delta\qq\nu\beta\gamma\ooo\delta\nu -
\frak S_{\alpha\beta\gamma}\qq\mu\alpha\beta\ooo 0\gamma .  \endalign
$$
 $\partial^2II$  also varies  as one moves along the fiber
(see [L2] to see how it varies).

\smallpagebreak

\noindent\S 9. {\bf Dimension of the tangential variety}

   The dimension of $\tau (X)$ is the dimension of its  tangent space at a
point.  Since we have adapted the frame bundle such that $\{A_0, A_{\alpha}\} =
\hat T_{[A_0]}$ we may take $[A_1]$ as a typical element of $\tau (X)$.  By the
defining equation of the Maurer-Cartan form,

$$
dA_1  \equiv  \ooo  {0} {1}A_0  +\ooo  {\rho} {1}A_{\rho} +
\ooo  {\mu} {1}A_ {\mu} \ \ \text{mod}\ \{ A_1 \} \tag 9.1
$$
and dim$T_{[A_1]}\tm = $ $\{$ the number of independent  1-forms in (9.1)
$\}$.  Over $X$,
$\{\ooo 0 1 ,\ooo\rho 1\}$ are independent so we only need to know the  number
of independent one-forms among the $\ooo{\mu}{1}$, but (up to twisting)
$$
\ooo  {\mu} {1}\ot (A_ {\mu}\text{mod}\hat T)\equiv  II_{\underline A_1}.
 \tag 9.2
$$
where $\underline A_1=\aa 0^*\ot (\aa 1\tmod \aa 0 ) \in T^*$.

   (9.2) implies
$$
\text{dim}\tm = n+ \text{dim}II_v(T) \tag 9.3
$$
  where $v$ is any
$II$-generic tangent
vector at a general point of $X$, proving (4.10).
Let $a_0:= \text {max}_{v\in T}(\text {dim}(II_v(T))$, so  dim$\tau (X) =
n+a_0$.

  Given $X^n\subset\pp {n+a}$, $\tm$ is degenerate if
dim(ker$II_v)>\text{max}\{ 0,n-a\}$ for generic (and therefore all) $v\in T$.
In other words, $\tm$ is degenerate
in the case $a\leq n$ if
for all $ [v]\in \Bbb P T$ there exists a $[Q]\in |II|$ such
that $[v]\in [Q]_{sing}$, and
$\tx$ is degenerate in the case  $a> n$
if there exists an $(a-n)$-dimensional
family of quadrics singular at $[v]$.
 By comparison, the Gauss map is degenerate if there exisits an
$[v]\in \Bbb PT$ such that every $[Q]\in |II|$ is singular at $[v]$.
 Of course $\tm$ is degenerate if the Gauss map of $X$ is degenerate, but  the
Gauss map of a smooth variety is nondegenerate unless $X$ is a linear space
(see [GH] or [L1]).

\smallpagebreak

\noindent\S 10. {\bf Dimension of the Secant Variety}

  Although by the connectedness theorem it is sufficient to work with the
tangential variety, we will gain additional information by studying the secant
variety as well. We now  infinitesimaly compute the dimension of  the secant
variety  $\sigma (X)$.  The computation is a straightforward generalization of
the computation in [GH] of the case where the codimension is large.  Let
 $[A_0]$ and $[B_0]$ denote points of $X$. Points of $\sigma (X)$ are of
the form
 $[A_0  + sB_0]$  where $s\in \Bbb C$.  To find the dimension of $\sigma (X)$
   we compute the number of linearly independent entries  in
 $$
 d(A_0 + sB_0) \equiv  \ooo  {0} {0}A_0 + \ooo  {\alpha} {}A_{\alpha} +
   s( \phi^0_0A_0 + \phi^{\alpha}_0A_{\alpha}) + dsB_0  \ \text{mod}(A_0 +
sB_0)
$$
where $\phi^A_B$ are the Maurer Cartan forms over $B_0$.

   We see the dimension of $\sigma (X)$ is one less than the number of linearly
independent vectors among

   $$A_0,A_{\alpha},B_0,B_{\beta}. $$
   (This gives an infinitesimal proof of Terracini's lemma (see [FL]) that if
$X^n\subset\pp N$ is smooth $x,y\in X$ and $p\in \overline{xy}-X$ are general
points then $ \hat T_p\sm = $Span$(\hat T_xX , \hat T_yX)$.)
   Now we will use the assumption that $X$ is smooth and connected.  Let
   $(B_0(t),\hdots , B_N(t))$ be a coframing over an arc $B_0(t)$ such that
$(B_0(0),\hdots , B_N(0)) = (A_0,\hdots , A_N) $ and $B_0(0)' = v$.      Let
$\Lambda_0 = \{A_0\hd  A_n\} $.
Expanding $B_0(t),B_{\alpha}(t)$ in Taylor series we have:
$$
\align
B_0(t) &= A_0 + tv + \frac{t^2}{2}\frp{{}^2B_0}{v\partial v}\mid_{t=0}
+ \frac{t^3}{3!}\frp{{}^3B_0}{v\partial v\partial v}\mid_{t=0} + O(t^4)\\
B_{\alpha}(t) &= A_{\alpha} + t\frp{B_{\alpha}}{v}\mid_{t=0}
 + \frac{t^2}{2}\frp{{}^2B_{\alpha}}{v\partial v}\mid_{t=0}
+ \frac{t^3}{3!}\frp{{}^3B_{\alpha}}{v\partial v\partial v}\mid_{t=0} +
O(t^4).\endalign
$$
Now mod$\Lambda_0$,
$$\align
\frp{{}^2B_0}{v\partial v}\mid_{t=0} &\equiv II(v,v)\\
\frp{B_{\alpha}}{v}\mid_{t=0} &\equiv II(v,\underline A_{\alpha})
\endalign
$$
and
mod$\{\Lambda_0, II_v(T)\}$,
$$
\align
\frp{{}^3B_{0}}{v\partial v\partial v}\mid_{t=0} &\equiv
III^v(v,v,v)\\
\frp{{}^2B_{\alpha}}{v\partial v}\mid_{t=0}&\equiv
III^v(v,v, \underline A_{\alpha}).\endalign
$$
(where  $\underline A_\alpha = \aa 0^*\ot(\aa\alpha\tmod\aa 0)\in T$).
Higher order terms in the series cannot contribute to the dimension of
$\sm$ because either the lower order terms (together with $\Lambda_0$) span
a space of maximal dimension and $\sm$ is nondegenerate, or $III^v\equiv 0$ and
since the higher fundamental forms lie in the prolongation of the lower
fundamental forms,
they are zero so the higher order terms cannot contribute any new directions.
So letting $v\in T$ be $II$-generic we obtain
$$
\text{dim}\sm = \left\{\matrix
n+ \text{dim}II_v(T) &\text{ if }III^v(v,v,v)= 0\\
 n+ \text{dim}II_v(T)  + 1 &\text{ if }III^v(v,v,v)\neq 0\endmatrix\right\}.
\tag 10.1
$$
(10.1) combined with (9.3) and  (1.3) implies:

\proclaim{Proposition 10.2} Let $X^n\subset\pp{n+a}$ be a smooth variety with
degenerate secant variety. Then for all $II$-generic vectors $v\in T$,
$III^v=0$ (and therefore  $III=0$ ).
Conversely, if $X^n\subset\pp{n+a}$ is a
variety and $III_X$ is not identically zero, then $\sx$ is nondegenerate.
\endproclaim

\proclaim{Theorem 10.3}
Let $Y\subset \ppp V$ be a variety and let
$X=v_d(Y)\subset\ppp S^dV$ be the Veronese re-embedding.
If $d>2$ or
$d=2$ and $Y$ is not a linear subspace of $\ppp V$,
then $\sx$ is nondegenerate.
\endproclaim

\demo{Proof}
By ([L2], 3.16), If
 $d>2$   or $d=2$ and  $Y$ is not a linear space, then
$III_{v_d(Y)}\neq 0$.
\enddemo

(10.3) extends a result of Roberts that follows from Proposition
2 of [R], where (10.3) is proved for $d>2$. It was also shown
by Kaji, [K], that (10.3) holds for homogeneous varieties.

One can describe the dimensions of higher secant varieties
in a similar manner
to (10.1) using fundamental forms of possibly higher order,
and make similar apriori statements like (10.3).   The
general   expression for $\sigma_k(X)$
is  complicated, but if
$III^v\equiv 0$ it becomes
$$
\text{dim}\sigma_k(X) = n+ \text{dim}\{ II_{v_1}(T)\hd II_{v_{k-1}}(T)\}\tag
10.4
$$
where $v_1\hd v_{k-1}$ are $II$-generic and in general position.
(The notation is such that $\sm = \sigma_2(X)$.)  In particular
 we recover the
following result of Zak and Fantechi [F]:

\proclaim{Theorem 10.5  }  Let $X^n\subset\pp{n+a}$ be a smooth variety
not contained in a hyperplane
with degenerate secant variety of dimension $n+a_0$.
Let $\sigma_k(X)$ denote the closure of the union of
all secant $\pp{k-1}$'s of $X$. Then
$$
\text{dim}\sigma_k(X)\leq n+ (k-1)a_0.
$$
\endproclaim

This estimate enables one to recover the other valid results of [Z2],
e.g. Theorem 3. Notice that to have $\tdim \sigma_k (X) = n+ (k-1)a_0$ is
rather difficult, so it might be interesting to study cases where equality
occurs.

\bigpagebreak

\noindent\S 11. {\bf Proofs of statements in \S 6}

Although to keep notation consistent we will continue to work on
frame bundles over $X$, all the results in this section are just
about systems of quadrics with defects. This is reflected in
 that we always work modulo  the semi-basic forms to prove the
statements in \S 6, which is the same as
saying
that we are working with a fixed system of quadrics.

We now assume $v=\underline A_1$ varies over $II$-generic vectors and adapt
frames to reflect the structure described in \S 6.
(Recall that we are using the notation
$\underline A_{\alpha} = A_0^*\ot (\aa\alpha\tmod\aa 0 )$.)
Introduce the index ranges
 $a_0+1= n- \tdim(\tker II_v ) +1
\leq \epsilon ,\delta \leq
n $, $a_0+1\leq h\leq a$ and restrict to frames
where
$$
\align \tann (v) &=
\{\ooo{n+h}\alpha\oo\alpha \}
\tag 11.1\\
\tker II_v &= \{\underline A_{\epsilon}\}. \tag 11.2\endalign
$$

  (11.1) implies
$$
\ooo{n+h} 1=0. \tag 11.3
$$
Introduce   additional index ranges $2\leq \rho ,\sigma \leq n$.
Using (8.10) (i.e. differentiating $\ooo {n+h} 1$), we have
$$
\rr{n+h}1\epsilon\beta\oo\beta = \qq{n+h}\epsilon\rho\ooo\rho 1 \tag 11.4
$$
Since  the forms
$\ooo\rho 1$
 are all linearly independent, and independent of the
semi-basic forms (since we are allowing $\aa 1$ to range over all
$II$-generic vectors), (11.4) implies
$$
 \qq {n+h}\epsilon\rho=0 \ \forall h,\epsilon ,\rho  \tag 11.5
$$
which reproves (6.5).
For later use we note that (11.4) also implies
$$
\rr {n+h} 1\epsilon\beta = 0 \ \forall h,\epsilon ,\beta .\tag 11.6
$$

\smallpagebreak

We reduce the frame bundle further by requiring
$$
\align \tsingloc (\tann ( v)) &= \{\underline A_1,\underline A_s ,\underline
A_{\epsilon}  \}
\tag 11.7 \\
T / \tsingloc (\tann ( v))  & =\{\underline A_j \}
\tmod \{\underline A_1,\underline A_s ,\underline
A_{\epsilon}  \} \tag 11.8  \endalign
$$
where $2\leq i,j,k,l\leq r+1$, $r+2\leq s,t\leq a_0$.

{}From now on, we will abuse notation, using $\aa B$ to denote
$\underline A_B$.

Let $1\leq \lambda \leq a-1$ so
$$
\{ \aa\lambda \} \tmod \{\aa\epsilon \}
 = \{ \aa 1,\aa j ,\aa s \}
\tmod \{\aa\epsilon \}
 = T/\tker (II_v)
.\tag 11.9
$$

Adapt frames to  the map $II_v=II_{ \aa 1}$ such that
$$
II_{ \aa 1}( \aa \lambda) = \aa{n+\lambda}\tag 11.10
$$

 On our bundle adapted to
(11.1),(11.2),(11.7),(11.8), and (11.10), which we will call $\cf^2$,
$$
\ooo{n+\lambda}1=\oo\lambda.\tag
11.11
$$

Our adaptations have the effect that
$$
\align
&\{ [Q]\in \ii\ | \ [v]\in [Q] \} = \ppp \{
\ooo{n+j}\alpha\oo\alpha ,
\ooo{n+s}\alpha\oo\alpha ,\ooo{n+h}\alpha\oo\alpha \}\tag 11.12\\
&
{F_v} = \ppp \{\aa{n+1}, \aa{n+s}
\tmod \hat T\}\tag 11.13\\
&
II^*({F_v}\upperp ) = \{
\ooo{n+j}\alpha\oo\alpha, \ooo{n+h}\alpha\oo\alpha \}\tag 11.14
\endalign
$$

\smallpagebreak

Continuing to use (8.10), i.e.  to take third derivatives, we  have
$$
\align
\rrn h11\alpha\oo\alpha &= -\ooo{n+h}{n+1}  \tag 11.15\\
\rrn h1j\alpha\oo\alpha &= -\ooo{n+h}{n+j} + \qq{n+h} ij\ooo i1 \\
\rrn h1s\alpha\oo\alpha &= -\ooo{n+h}{n+s}  \endalign
$$

Introduce further index ranges
$
\xi ,\eta \in \{ 1,s,\epsilon \}
$,
i.e., let  $\{ A_{\xi} \}$ be such that
$$
\{\underline A_{\xi} \} = \tsingloc (\tann (v)). \tag 11.16
$$

\demo{Proof of (6.16)}
Using (8.10) again
$$
\rrn h \xi\eta\beta\oo\beta = - \qqn \lambda\xi\eta\ooo{n+h}{n+\lambda}.
\tag 11.17
$$

Substituting (11.15) into (11.17) implies
$$
0\equiv - \qqn j\xi\eta\qq{n+h} ji\ooo i 1 \tmod \{\oo\beta \}
\ \forall h. \tag 11.18
$$
Since $\tsingloc ( \{ \qq{n+h}ij\oo i\oo j \}|_{\{\aa k \} })=0$,
(11.18)
  implies
$$
\qqn j\xi\eta =0\  \forall j,\xi,\eta \tag 11.19
$$
Comparing (11.14),(11.16) and (11.19) proves (6.16).\qed
\enddemo

\demo{Proof of (6.21)}
Restrict to the case $\tx$ is a hypersurface
 (so $a_0=a-1, h=a$ and $1\leq\lambda\leq a-1$). Continuing to use
(8.10),

 $$
\align &
\rrn a k\epsilon\beta\oo \beta = -\qqn jk\epsilon\ooo{n+a}{n+k}
     -\qq{n+t}k\epsilon\ooo{n+a}{n+t}+\ooo k\epsilon \tag 11.20
\\
&
\rrn j1\epsilon\beta\oo\beta = \ooo j\epsilon + \qqn j\epsilon k \ooo k1
 \tag 11.21\endalign
$$

Combining (11.20), (11.21) with (11.15) yields
$$
\align
\ooo j\epsilon&\equiv \qqn jk\epsilon\ooo k1 \tmod \{\oo\alpha \} \tag 11.22\\
\qq\na ij\ooo j\epsilon
 &\equiv -\qqn k i\epsilon \ooo\na{n+k}\tmod \{\oo\alpha
\}  \\
 &
\equiv -\qqn k i\epsilon  (\qq\na kl \ooo l1 )
\tmod \{\oo\alpha \} \tag 11.23\endalign
$$
combining (11.22),(11.23) yields
$$
\qq\na ij\qqn j i\epsilon + \qqn k i\epsilon\qq\na kj =0.\tag 11.24
$$
Here $P_v = \ooo\na\alpha\oo\alpha$. We   restrict to frames where
$$
\qq\na ij = \delta_{ij} \tag 11.25
$$
and (11.24) becomes
$$
 \qqn j i\epsilon + \qqn ij\epsilon  =0 \tag 11.26
$$
proving (6.21).\qed
\enddemo

\demo{Proof of (6.23)}
By the assumption that $Z$ is not built out of a tangent
developable,
we can make a further adaptation that
$\ppp \{\aa{n+s}\}=V_Z$, i.e. that
$\{\aa{n+s}\}\tmod \hat T$
is a fixed subspace as
$v$ varies in $T^{gen}$.
We need to work with the subsystem $II^*(V_Z\upperp )$.
The notationally easiest way to do this is just set $\{\aa{n+s}\}
=\{ 0\}$ in everything above. There is a potential
for confusion here, because for the following argument only,
$r+2\leq\epsilon ,\delta\leq n$, i.e.
$\{\aa 1 ,\aa\epsilon
\}$ will denote
$\tsingloc (\tann (v))$
where   $\{\aa 1 ,\aa s ,\aa\epsilon \}$ denoted
that same space above. On the smaller
system $II^*(V_Z\upperp )$, the spaces $\tsingloc(\tann (v))$ and
$\{ v,\tker II_v \}$ coincide.

\smallpagebreak

Again, using (8.10),
$$
\rrn j 11\beta\oo\beta = -\ooo{n+j}{n+1}
+ 2\ooo j1 \tag 11.27
$$
i.e.
$$
\ooo{n+j}{n+1}
\equiv  2\ooo j1 \tmod \{\oo\beta \}. \tag 11.28
$$
Finally
$$
\align
\rrn j \epsilon\delta\beta\oo\beta &=
-\qqn 1\epsilon\delta\ooo{n+j}{n+1}+ \qqn j\epsilon k \ooo k\delta + \qqn
j\delta k \ooo k\epsilon \\
 &\equiv
- (2\qqn 1\epsilon\delta\delta_{ij} + \qqn j\epsilon k \qqn k i\delta + \qqn
j\delta k\qqn k i\epsilon)\ooo i1 \tmod \{ \oo\alpha \} \tag 11.29
\endalign
$$
(where we have used (11.22) and (11.28)). (11.29) implies
$$
\qqn j\epsilon k\qqn k\delta i + \qqn j\delta k\qqn k\epsilon i
+ 2\qqn 1\epsilon\delta \delta^j_i =0
\ \forall \epsilon,\delta, i,j \tag 11.30
$$
proving (6.23).\qed\enddemo

\bigpagebreak

\noindent\S 12. {\bf A refinement of a result in [GH] }

 \proclaim {Theorem 12.1} Let $X^n\subset\pp {n+a}$
be a
 variety (respectively a patch of a complex manifold)
 with degenerate tangential variety (resp. manifold) of
dimension
$n+a_0$.  Then
 the Gauss map of $\tx$ has fibers of dimension at least
$\delta_{\tau}(X)+1$, where $\delta_{\tau}(X)$
is the tangential defect of $X$.
 \endproclaim

(12.1) refines ([GH], 5.17) where it was observed that
 the Gauss map of
a degenerate
$\tx$ has fibers of dimension at least
$2$.

\demo{Proof}
The dimension of a fiber of $(\gamma (\tau (X_{sm})))$ is
 $\tdim(\tsingloc |II_{\tm [\aa 1]}|)$.

{}From (9.1), we   compute
the second
fundamental form of $\tx$ at $[\aa 1]$:
$$
II_{\tm [\aa 1]}
= (\ooo{n+h}j\ooo j1 + \ooo{n+h}{n+\lambda}\ooo{n+\lambda}1)
\ot\aa{n+h}\tmod\{\aa 0\hd\aa{n+a-1}\}.
\tag 12.2
$$
However
(11.15) combined with (11.6) implies
$\ooo{n+h}{n+\lambda}\equiv 0\tmod \{ \oo\lambda ,\ooo j1 \}$.

Thus $\tdim(\tsingloc |II_{\tm [\aa 1]}|)
\geq (n+a_0)-\tdim\{\oo\lambda ,\ooo j1\}
=(n+a_0)-(a_0-1+r)=n-r+1$.
But now $r\leq a_0-1$ and $\delta_{\tau}(X)=n-a_0$.
 \qed\enddemo

   Any tangential variety is ruled by parallel lines.    (12.1) says that if a
tangential variety is
degenerate, it is ruled by parallel
$(n-a_0+1)$-planes. (12.1) will be strengthened in the case $X$ is
a
smooth variety.

\bigpagebreak

\noindent \S 13. {\bf Global Results}

We will now restrict to the case $X^n\subset\pp{n+a}$ is a smooth variety with
degenerate tangential variety that is a hypersurface, and
use the    global  information that $X$ is smooth to assume
$III^v=0$ via (10.2). By (8.7),
$$
\align
III^v &=
\ooo{n+a}{n+\lambda}\ooo{n+\lambda}\xi\oo\xi\mid_{S^3(SA_1)}\ot\aa{n+a}\tmod
\{\aa 0\hd\aa{n+a-1}\}
\tag 13.1\\
&= -\rrn a\xi\eta\zeta\oo\xi\oo\eta\oo\zeta
\ot\aa{n+a}\tmod
\{\aa 0\hd\aa{n+a-1}\}\endalign
$$
where in the second line we have used (11.17). So
$III^v=0$ implies
$$
  \rrn a\xi\eta\zeta =0\ \forall \xi,\eta ,\zeta .\tag 13.2
$$

\smallpagebreak

\demo{Proof of (7.9)}
 Let $A^B$ denote the dual basis
element to $A_B$ at a point $f=(\aa 0\hd \aa\na )$ of $\cf^2$.
Consider
$$
\align \tilde\rho : \cf^2&\ra X^*_{\Delta}\tag 13.3\\
(\aa 0\hd\aa{n+a})&\mapsto [ A^{n+a}].\endalign
$$
The only restriction on $[ A^{n+a}]$ comes from (11.1), so
$\tilde\rho$ is surjective. Note that
$$
\tilde\rho _{* f} = -\ooo \na j
\ot A^j
-
\ooo \na {n+\lambda}
\ot A^{n+\lambda}
\tmod
\{ A^{\na}\} .\tag 13.4
$$
By (13.4), the rank of $\tilde\rho _{* f}$ equals the number of
independent forms among the $\{ \ooo\na j ,\ooo\na{n+\lambda}\}$.

 Using (8.10)
further we have
$$
 \align
\rrn a11j \oo j &= -\ooo{n+a}{n+1} \tag 13.5\\
\rrn a 1sj\oo j &= -\ooo{n+a}{n+s}\\
\rrn a 1 k\lambda\oo\lambda &= -\ooo{n+a}{n+k} +\ooo k 1.
\endalign
$$
where we have used (13.2)   and (11.6).

To facilitate the computations, we will restrict to frames where
some of the coefficients of $\pii$ are zero. In fact we use
the freedom in
$\gg k{n+1} , \gg k{n+s}$ (see (8.14)) to restrict to frames where
$$
\rrn a 11k , \rrn a 1sk =0. \tag 13.6
$$
Call the new bundle $\cf^3$. On $\cf^3$,
$$
\{\ooo\na{n+\lambda} ,\ooo\na j \} = \{ \ooo j1 ,\oo j \}
\tag 13.7
$$

Since  the projections $\cf^3\ra X$, $\cf^3\ra X^*_{\Delta}$
factor through $\ci_{\Delta}$, we have
$$
\trank (\tilde\rho|_{\cf^3})_{* f} = \trank (\rho_{\Delta})_{* ([\aa 0],
[A^{n+a}])}. \tag 13.8
$$
But by  (13.7), $\trank (\tilde\rho|_{\cf^3})_{* f} = 2r$.\qed\enddemo

\smallpagebreak

The computation above also allows us to improve (12.1) under
the assumption that $X$ is a smooth variety.

Using the computations above, (12.2)   reduces to
$$
\align
II_{\tm \, [\aa 1]} &=
(\ooo{n+a}j\ooo j1+ \ooo{n+a}{n+\lambda}\ooo{n+\lambda}1)\ot
(\aa{n+a}
\tmod \{\aa
0\hd\aa
{n+a-1}\})\tag 13.9\\
&= (2\delta_{jk}-\rr\na 1jk )\ooo j 1\oo k\ot  (\aa{n+a}\tmod \{\aa 0\hd\aa
{n+a-1}\}). \endalign
$$

 \proclaim {Theorem 13.10} Let $X^n\subset\pp {n+a}$
be a smooth variety with degenerate tangential variety of dimension $n+a_0$.
Then
the Gauss map of $\tx$ has
fibers of dimension at least $\delta_{\tau}(X)+2$,
where $\delta_{\tau}(X)$ is the tangential defect of $X$.
 \endproclaim

\demo{Proof} If $\tx$ is not
a hypersurface, we can project $X$   smoothly to a smaller projective
space where it is, and the dimension of the Gauss image of
$\tx$ will not change, and (13.9) proves
the case $\tx$ is a hypersurface.
\enddemo

  (13.10) is sharp.  Equality holds for Severi varieties  and
  that strict inequality occurs for the Grassmannian
$G(2,7)\subset\pp {20}$.

\bigpagebreak

\noindent\S 14.  {\bf Zak's Theorem on Severi Varieties}

\smallpagebreak

Assume $a=\frac n2+2$. (7.3) implies
$\gamma_Z$ is nondegenerate at general points.

Assuming $\aa 2$ is $II$-generic (which is allowed since $\aa 2\in \{\aa j\}$),
use
$\gg{n+1}{n+k},\gg{n+1}\na$ to restrict to frames where
$$
[q^{n+1}]=[\qq{n+1}\alpha\beta\oo\alpha\oo\beta ]= \tann (\aa 2 ).
\tag 14.1
$$
 This implies
$\trank ( q^{n+1} )=r$.
Using the Clifford algebra structure and the rank restriction theorem,
one can check that $q^{n+1}|_{\{\aa 1,\aa\epsilon\} }$ is non-degenerate.
(The $n=2, 4$ and $6$ cases are clear, for the $n=8$ case one would
have a non-zero $4\times 4$ skew symmetric matrix that squared to zero
(which is possible) but that also skew commuted with
a different non-zero
$4\times 4$ skew symmetric matrix that either squared to be the identity
or zero
(which is not possible). One argues similarly for larger $n$.)
After
using the remaining freedom in $\gg\epsilon j$, we may (and do)
restrict to frames where
$$
\tsingloc (q^{n+1}) =\{\aa j \} \tag 14.2
$$
(Note that (14.2) makes sense
since the $\gg 1k$ are
tied to the $\gg{n+1}{n+a}$, so
 $\{\aa j \}$ is
now  a well defined subspace of $T$.)

Because of the restrictions on dimensions of Clifford modules,
$n=2,4,8$ or $16$. (In the
$n=2$ case,  $\{\aa j \}=\{\aa 2\}$
and $\{\aa\epsilon \}=\{ 0\}$.)

   Using some of the freedom in
$\gg\epsilon\delta  $ we may restrict to frames where
$$
\qqn 1\epsilon\delta = \delta_{\epsilon\delta}.
  \tag 14.3
$$

By using our  remaining freedom in $g^k_j$ and $g^{\epsilon}_{\delta}$
we may arrange that   $\phi (A_{\epsilon}) = \qqn i j\epsilon\oo j\otimes A_i$
corresponds to $\Bbb A$-multiplication by
$J_{\epsilon -(\frac{n}{2} +1)}$ on the vector space $\{
A_j\}$, where the $J$'s are as in \S 3.  We are
guaranteed that we are able to do this because  the Clifford representation is
nondegenerate, and in these dimensions there are unique irreducible
representations.

  (6.16) applied to $\aa 2$ implies   that for all $i,j,k$,
$$
\qqn ijk =0. \tag 14.4
$$
Let $J_{\epsilon}(j)$ be the unique index $k$ such that
$\qqn j{\epsilon}k = \pm 1$ and we use the sign convention that if
$\qqn j {\epsilon} k = -1$ then we take $\oo {J_{\epsilon}(j)} = -\oo k$.
We have
$$
II= \pp{}\{ (\oo 1)^2 + \Sigma_{\epsilon}(\oo {\epsilon})^2, \oo 1 \oo j +
\Sigma_{\epsilon}\oo {\epsilon}\oo {J_{\epsilon}(j)}, \Sigma_j(\oo j)^2\} \tag
14.5
$$
 which is seen to be
the second fundamental form of a Severi variety by comparing with (5.2).
We also have
$$
0= |III^{\aa j}| =| (\ooo{n+1}{n+k}\ooo{n+k}i\oo i +
\ooo{n+1}{n+a}\ooo{n+a} j\oo j)\mid_{\tsingloc (\tann (A_j))}| =
|\rrn 1jki\oo i\oo j\oo k | \tag 14.6
$$
giving
$$
\rrn 1 jki=0.\tag 14.7
$$
Restrict frames futher to normalize away more of $\pii$ by
using $\gg 1{n+1},\gg 1{n+k}, \gg \epsilon{n+1},\gg \epsilon{n+k}$
respectively  to set
$$
\rrn 1111 = 0,\  \rrn 11k1 = 0,
\rrn 111\epsilon = 0,\  \rrn 11k\epsilon = 0. \tag 14.8
$$

Fixing any $k$, normalize away $\rrn k1kk $ (and therefore
$\rrn j1jj\ \forall  j$) by using
$
\gg k\na
$.

In summary, combining all the adaptations on our new frame bundle;
$$
\align
 0
&= -\ooo 0 0 - \ooo{n+1}{n+1} + 2\ooo 1 1\tag 14.9\\
\rrn 1 1 j\epsilon\oo\epsilon
&= - \ooo{n+1}{n+j} + \ooo 1 j \\
\rrn 11\epsilon\delta\oo\delta
&= \ooo\epsilon 1 + \ooo 1\epsilon \\
\rrn 1jk\epsilon\oo\epsilon
&= - \delta_{jk}\ooo{n+1}{n+a}\\
\rrn 1 j\epsilon\rho\oo\rho
&= - \ooo{n+1}{n+J_{\epsilon}(j) }+ \ooo\epsilon j \\
\rrn 1 \epsilon\delta\beta\oo\beta
&
= -\delta_{\epsilon\delta}(\ooo 0 0 + \ooo{n+1}{n+1}) + \ooo\epsilon\delta
+\ooo\delta\epsilon\\
\rrn k 11\beta\oo\beta
&= -\ooo{n+k}{n+1} + 2\ooo k 1\\
\rrn k1 j\beta\oo\beta
&=
\delta^j_k (\ooo 1 1 - \ooo 0 0 ) -
\ooo{n+k}{n+j} + \ooo kj + \qqn kij\ooo i1  +
\delta^{J_{\epsilon}(k)}_j\ooo\epsilon 1\\
\rrn k 1\epsilon\beta\oo\beta
&= \ooo k\epsilon + \ooo{J_{\epsilon}(k)} 1\\
\rrn k ij\beta\oo\beta &= -\delta_{ij}\ooo{n+k}{n+a} + \qqn
ki\epsilon\ooo\epsilon j
+\qqn kj\epsilon\ooo\epsilon i\\
\rrn k j\epsilon\beta\oo\beta
&= -\delta^{J_{\epsilon}(k)}_j\ooo 00 +
\ooo {n+k}{n+J_{\epsilon} (j)} + \ooo{J_{\epsilon}(k) }j +
\delta^{J_{\delta}(k)}_j\ooo\delta\epsilon\\
\rrn k\epsilon\delta\beta\oo\beta
&= -\delta_{\epsilon\delta}\ooo{n+k}{n+1}
+\ooo{J_{\epsilon} (k)}\delta + \ooo{J_{\delta}(k)}\epsilon\\
\rrn ajk\beta\oo\beta
&= -\delta_{jk}(\ooo 0 0 + \ooo{n+a}{n+a} ) +
\ooo jk + \ooo kj\\
&0=\ooo{n+a}{n+1}\\
&0=\ooo{J_{\epsilon}(k)}1 + \ooo k\epsilon .
\endalign
$$

Examing the terms in (14.9) combined with our new adaptations, we find that
all the terms $\rr\mu\alpha\beta\gamma=0$ and the   relations
imposed on the $\{ \ooo\alpha\beta ,\ooo\mu\nu \} $ are the ones
obtained by setting the right hand side of (14.9) to zero.
(Perhaps a better way to see all the coefficients of $\pii$ must
be zero is,
lettting $s=(s^1\hd s^n)$,   use the fact that
$III^{s^{\alpha}\aa\alpha}\equiv 0$
for generic and therefore all choices of $s$. One gets a polynomial in the
$s^{\alpha}$ with coefficients the $\rr\mu\alpha\beta\gamma$ that must
be zero, which together with our
normalizations, forces all coefficients of $\pii$ to be zero.)
Examining the $\partial^2II$ equations, one gets that the
$\rr\mu\alpha\beta{\gamma\delta}$ are normalizable to zero by motions of the
form $\aa\mu\ra\aa\mu + g^0_{\mu}\aa 0$ and one obtains relations giving the
forms $\ooo\alpha\mu$ in terms of the $\ooo 0\beta$.  Finally,  the equation
for the
coefficients of $\partial^3II$ reduce to
$\rr\mu\alpha\beta{\gamma\delta\epsilon}\oo\epsilon =
\frak S_{\alpha\beta\gamma\delta}\qq\mu\alpha\beta\qq\nu\gamma\delta
\ooo 0\nu$
which (using our explicit knowledge of $II$) tells us $\ooo 0\nu =0$
and $\rr\mu\alpha\beta{\gamma\delta\epsilon}=0$ .  One is left with
the structure equations for the corresponding Lie group proving the variety is
indeed the
expected homogeneous space.
(One gets an aesthetically nicer basis of the Lie algebras if one switches
the numberings of the $\aa\epsilon$ and the $\aa j$.)\qed

\medpagebreak

{\bf References}:

\smallpagebreak

\noindent [D] Degoli, L., {\it Due Nouvi Teoremi Sui Sistemi Lineari
Di Quadriche A Jacobiana Identicamente Nulla},
Collect.-Math 35(1984), pp. 125-139.

\smallpagebreak

\noindent [ESB] Ein, L. and N. Shepherd-Barron, {\it Some Special Cremona
Transformations}, Amer. J. Math {\bf 111} (1989), pp.783-800.

\smallpagebreak

\noindent [F] Fantechi, B., {\it On the Superaddivity of Secant Defects},
Bull.
Soc. Math. France, {\bf 118}(1990) p. 85-100.

\smallpagebreak

\noindent [FL] Fulton, W. and R. Lazarsfeld, {\it Connectivity and its
applications in algebraic geometry}, in Algebraic Geometry (Proceedings). Lect.
Notes in Math., No. {\bf 862}, pp. 26-92.

\smallpagebreak

\noindent [GH] Griffiths, P. and J. Harris, {\it Algebraic Geometry and Local
Differential Geometry}, Ann. scient. Ec. Norm. Sup.  {\bf 12} (1979), pp.
355-432.

\smallpagebreak

\noindent [Ht] Hartshorne, R., {\it Varieties of small codimension in
projective space}, Bull. A.M.S. {\bf 80} (1974), pp. 1017-1032.

\smallpagebreak

\noindent [Hv] Harvey, F. R., {\it Spinors and Calibrations}, Academic Press
(1990).

\smallpagebreak

\noindent [K] Kaji, H. {\it On the homogeneous
projective varieties with degenerate secants}. preprint.

\smallpagebreak

\noindent [L1] Landsberg, J.M., {\it On Second Fundamental Forms
of Projective Varieties}.  Invent. Math.  {\bf 117}(1994) pp. 303-315.

\smallpagebreak

\noindent [L2] Landsberg, J.M.,
{\it  Differential-geometric characterizations of
complete intersections. }
	 	 To appear in {   Journal of Differential Geometry}. (alg-geom/9407002
eprint)

\smallpagebreak

\noindent [LM] Lawson, H.B. and
M. Michelsohn,  {\it Spin Geometry}, Princeton
University Press (1989).

\smallpagebreak

\noindent [LV]  Lazarsfeld, R. and A. Van de Ven {\it Topics in the Geometry of
Projective Space, Recent Work of F.L. Zak}, DMV Seminar (1984), Birkhauser.

\smallpagebreak

\noindent [R]  Roberts, J. {\it Generic projections of algebraic
varieties}. Amer. J. Math. , {\bf 93}(1971) pp. 191-214.

\smallpagebreak

\noindent [T]  Terracini, {\it Alcune questioni sugli spazi tangenti e
osculatori ad una varieta, I, II, III}. Atti
Della Societa dei
Naturalisti e Matematici, (1913) pp. 214-247.

\smallpagebreak

\noindent [Z1]  Zak, F.L., {\it Tangents and Secants of Algebraic Varieties},
AMS Translations of mathematical monographs, {\bf 127}(1993).

\smallpagebreak

\noindent [Z2]  Zak, F.L., {\it Linear Systems of Hyperplane sections of
varieties of low codimension},  Functional Anal. Appl.,
t. {\bf 19}(1985), pp. 165-173.

\enddocument